\documentclass[12pt,twoside,a4paper]{article}
\usepackage{bm}
\usepackage{braket}
\usepackage{graphicx}
\usepackage{graphics}
\usepackage{theorem}
\usepackage{amsmath}
\usepackage{amssymb,mathrsfs}
\usepackage{latexsym}
\usepackage{cases}
\usepackage[dvipdfmx,bookmarksnumbered,colorlinks,linkcolor=blue,urlcolor=blue,citecolor=blue]{hyperref}
\usepackage{color}
\usepackage{algorithm}
\usepackage{algpseudocode}
\usepackage{caption}
\newfloat{procedure}{H}{loa}
\floatname{procedure}{Procedure} 
\newtheorem{THM}{Theorem}[section]
\newtheorem{LEM}[THM]{Lemma}
\newtheorem{PROP}[THM]{Proposition}
\newtheorem{COR}[THM]{Corollary}

\newtheorem{ASSU}[THM]{Assumption}

\newtheorem{REM}[THM]{Remark}
\newtheorem{DEF}[THM]{Definition}

\newtheorem{EXA}[THM]{Example}

\newtheorem{hm-cond}{Condition}
\usepackage[dvipdfmx]{hyperref}
\hypersetup{
bookmarksnumbered,colorlinks=false,linkcolor=blue,urlcolor=blue,citecolor=blue
}

\numberwithin{equation}{section}  
\newcommand{\vc}{\bm}

\makeatletter
\DeclareRobustCommand\widecheck[1]{{\mathpalette\@widecheck{#1}}}
\def\@widecheck#1#2{%
    \setbox\z@\hbox{\m@th$#1#2$}%
    \setbox\tw@\hbox{\m@th$#1%
       \widehat{%
          \vrule\@width\z@\@height\ht\z@
          \vrule\@height\z@\@width\wd\z@}$}%
    \dp\tw@-\ht\z@
    \@tempdima\ht\z@ \advance\@tempdima2\ht\tw@ \divide\@tempdima\thr@@
    \setbox\tw@\hbox{%
       \raise\@tempdima\hbox{\scalebox{1}[-1]{\lower\@tempdima\box
\tw@}}}%
    {\ooalign{\box\tw@ \cr \box\z@}}}
\makeatother
\makeatletter
\newif\if@borderstar
\def\bordermatrix{\@ifnextchar*{%
 \@borderstartrue\@bordermatrix@i}{\@borderstarfalse\@bordermatrix@i*}%
}
\def\@bordermatrix@i*{\@ifnextchar[{\@bordermatrix@ii}{\@bordermatrix@ii[()]}}
\def\@bordermatrix@ii[#1]#2{%
\begingroup
 \m@th\@tempdima8.75\p@\setbox\z@\vbox{%
 \def\cr{\crcr\noalign{\kern 2\p@\global\let\cr\endline }}%
 \ialign {$##$\hfil\kern 2\p@\kern\@tempdima & \thinspace %
  \hfil $##$\hfil && \quad\hfil $##$\hfil\crcr\omit\strut %
  \hfil\crcr\noalign{\kern -\baselineskip}#2\crcr\omit %
  \strut\cr}}%
 \setbox\tw@\vbox{\unvcopy\z@\global\setbox\@ne\lastbox}%
 \setbox\tw@\hbox{\unhbox\@ne\unskip\global\setbox\@ne\lastbox}%
 \setbox\tw@\hbox{%
  $\kern\wd\@ne\kern -\@tempdima\left\@firstoftwo#1%
  \if@borderstar\kern 2pt\else\kern -\wd\@ne\fi%
 \global\setbox\@ne\vbox{\box\@ne\if@borderstar\else\kern 2\p@\fi}%
 \vcenter{\if@borderstar\else\kern -\ht\@ne\fi%
  \unvbox\z@\kern -\if@borderstar2\fi\baselineskip}%
 \if@borderstar\kern-2\@tempdima\kern2\p@\else\,\fi\right\@secondoftwo#1 $%
 }\null \;\vbox{\kern\ht\@ne\box\tw@}%
\endgroup
}
\makeatother

\newcommand{\ol}{\overline}
\newcommand{\ool}[1]{\overline{\overline{\bm{#1}}}}

\newcommand{\wt}{\widetilde}

\newcommand{\down}[2]{\smash{\lower#1\hbox{#2}}}
\newcommand{\up}[2]{\smash{\lower-#1\hbox{#2}}}

\newcommand{\dm}{\displaystyle}
\newcommand{\qed}{\hspace*{\fill}$\Box$}
\newcommand{\proof}{\noindent {\it Proof:~}}

\newcommand{\vmin}{\wedge}

\newcommand{\EE}{\mathbb{E}}
\newcommand{\PP}{\mathbb{P}}

\newcommand{\calL}{\mathcal{L}}

\newcommand{\calS}{\mathcal{S}}

\newcommand{\bcal}[1]{\bm{\mathcal{#1}}}

\newcommand{\bbL}{\mathbb{L}}
\newcommand{\bbM}{\mathbb{M}}
\newcommand{\bbN}{\mathbb{N}}
\newcommand{\bbR}{\mathbb{R}}
\newcommand{\bbS}{\mathbb{S}}

\newcommand{\bbZ}{\mathbb{Z}}

\newcommand{\re}{{\rm e}}

\newcommand{\dd}[1]{\if#11 1\!\!1
\else {\if#1C I\!\!\!C
\else {\if#1G I\!\!\!G
\else {\if#1J J\!\!\!J
\else {\if#1S S\!\!\!S
\else {\if#1Z Z\!\!\!Z
\else {\if#1Q O\!\!\!\!Q
\else I\!\!#1
\fi}
\fi}
\fi}
\fi}
\fi}
\fi}
\fi}
\makeatletter
\def\widebar{\accentset{{\cc@style\underline{\mskip10mu}}}}
\def\Widebar{\accentset{{\cc@style\underline{\mskip8mu}}}}
\makeatother
\newcommand{\wbartil}[1]{\if#1L \widebar{\hspace{-0.12zw}\widetilde{#1}}
\else {\if#1M \widebar{\widetilde{\!M\!}}
\else {\if#1W \widebar{\widetilde{\!W\!}}
\else {\if#1U \widebar{\hspace{-0.03zw}\widetilde{#1}}
\else {\if#1V \widebar{\hspace{-0.03zw}\widetilde{#1}}
\else {\if#1Y \widebar{\hspace{-0.0zw}\widetilde{#1}}
\else \,\widebar{\!\widetilde{#1}}
\fi}
\fi}
\fi}
\fi}
\fi}
\fi}
\makeatletter
\def\eqnarray{\stepcounter{equation}\let\@currentlabel=\theequation
\global\@eqnswtrue
\global\@eqcnt\z@\tabskip\@centering\let\\=\@eqncr
$$\halign to \displaywidth\bgroup\@eqnsel\hskip\@centering
$\displaystyle\tabskip\z@{##}$&\global\@eqcnt\@ne
\hfil$\;{##}\;$\hfil
&\global\@eqcnt\tw@ $\displaystyle\tabskip\z@{##}$\hfil
\tabskip\@centering&\llap{##}\tabskip\z@\cr}
\makeatother

\makeatother

\begin{document}\thispagestyle{empty}

\hfill

\vspace{-10mm}

{\large{\bf
\begin{center}
LEVEL-WISE SUBGEOMETRIC CONVERGENCE OF
THE LEVEL-INCREMENT TRUNCATION APPROXIMATION OF
M/G/1-TYPE MARKOV CHAINS%
\footnote[1]{%
To appear in Advances in Journal of the Operations Research Society of Japan, vol. 65, no. 4, October 2022
}
%
%
\end{center}
}
}

\begin{center}
{
\begin{tabular}[h]{cc}
Katsuhisa Ouchi\footnotemark[2] & Hiroyuki Masuyama\footnotemark[3]                        \\ 
\textit{Kyoto University}&\textit{Tokyo Metropolitan University}\\ 
\end{tabular}
\footnotetext[2]{E-mail: o-uchi@sys.i.kyoto-u.ac.jp}
\footnotetext[3]{E-mail: masuyama@tmu.ac.jp}
}

\bigskip
\medskip

{\small
\textbf{Abstract}

\medskip

\begin{tabular}{p{0.85\textwidth}}
This paper considers the level-increment (LI) truncation approximation of M/G/1-type Markov chains. The LI truncation approximation is useful for implementing the M/G/1 paradigm, which is the framework for computing the stationary distribution of M/G/1-type Markov chains. The main result of this paper is a subgeometric convergence formula for the total variation distance between the original stationary distribution and its LI truncation approximation. Suppose that the equilibrium level-increment distribution is subexponential, and that the downward transition matrix is rank one. We then show that the convergence rate of the total variation error of the LI truncation approximation is equal to that of the tail of the equilibrium level-increment distribution and that of the tail of the original stationary distribution.
\end{tabular}
}
\end{center}

\begin{center}
\begin{tabular}{p{0.90\textwidth}}
{\small

\medskip

}
\end{tabular}

\end{center}

\section{Introduction}
\label{sec:Intro}
The M/G/1 paradigm \cite{Neut89} is the framework for computing the stationary distribution in M/G/1-type Markov chains. The classical theory on M/G/1-type Markov chains \cite{Neut89} is useful for the algorithmic analysis of various semi-Markovian queues such as BMAP/GI/1 queues (see, e.g., \cite{Luca91,Neut89,Taki00}) and their multiclass extensions (see, e.g., \cite{Masu03-STM,Taki01-QUESTA37,Taki01-QUESTA39}). Such semi-Markovian queues have been comprehensively studied by many researchers over thirty years.

Implementing the M/G/1 paradigm needs the {\it level-increment (LI) truncation approximation}. The LI truncation approximation transforms the original transition probability matrix into another M/G/1-type stochastic matrix, which specifies an (M/G/1-type) Markov chain having level increments truncated at some upper bound. The resulting M/G/1-type stochastic matrix (resp.\ its corresponding M/G/1-type Markov chain) is referred to as a {\it level-increment (LI) truncation approximation} to the original M/G/1-type transition probability matrix (resp.\ the original M/G/1-type Markov chain).

The stationary distribution of the LI truncation approximation to the original chain can be computed by the M/G/1 paradigm, as an approximation to that of the original chain. This is because the LI truncation approximation is characterized with a finite number of block component matrices. Hence, we use the name ``an LI truncation approximation to the original stationary distribution" for the stationary distribution of an LI truncation approximation of the M/G/1-type Markov chain.

To the best of our knowledge, there are no previous studies on the evaluation of the error of an LI truncation approximation to the original stationary distribution. Thus, we have had no idea how the truncation parameter (the maximum of upward jumps in level) should be set in advance to meet a given error tolerance. In other words, we have to rely on our empirical knowledge to set the truncation parameter and naturally we cannot guarantee the accuracy of LI truncation approximations to be computed.

This paper studies the case where the error of LI truncation approximations to the original stationary distribution of the M/G/1-type Markov chain decreases at a subgeometric rate, that is, at a rate much slower than exponential, such as $N^{-\alpha}$ and $\exp\{-N^{\beta}\}$ ($N$ is the truncation parameter). In such a case, it takes tremendous computational cost to meet the error tolerance. Therefore, it is important to identify in what cases such an unfavorable situation occur and to what extent its impact is.

In fact, there is a study \cite{Masu21-M/G/1-Subexp} on the above-mentioned unfavorable cases for the {\it last-column-block-augmented (LCBA) truncation approximation} of M/G/1-type Markov chains. The LCBA truncation approximation to an M/G/1-type Markov chain can be considered a finite-level M/G/1-type Markov chain (see \cite[Remark~2.2]{Masu21-M/G/1-Subexp}). The study \cite{Masu21-M/G/1-Subexp} presents a subgeometric convergence formula for the level-wise difference between the original stationary distribution and its LCBA truncation approximation. The subgeometric convergence formula requires two conditions: (i) the finiteness of the second moment of level increments; and (ii)  the subexponentiality of the integrated tail distribution of nonnegative level increments in steady state. Besides, there are several studies \cite{Masu15-AAP,Masu16-SIMAX,Masu17-LAA,Masu17-JORSJ} on upper bounds for the error of LCBA truncation approximations to block structure Markov chains including M/G/1-type Markov chains.

The main contribution of this paper is to derive a subgeometric convergence formula for the {\it level-wise (not whole)} difference between the original stationary distribution and its LI truncation approximation. Our subgeometric convergence formula is in the same form as the one presented in \cite{Masu21-M/G/1-Subexp}, but our assumption is weaker than the above-mentioned conditions (i) and (ii) assumed in \cite{Masu21-M/G/1-Subexp}. We note that, as with the one in \cite{Masu21-M/G/1-Subexp}, our subgeometric convergence formula shows that the total variation norm of the relative {\it level-wise} difference between the original stationary distribution and its LI truncation approximation is asymptotically independent of the level variable. Moreover, this fact implies that there are no bias (in an asymptotic sense) in the convergence speeds of the LI truncation approximations to the level-wise subvectors of the stationary distribution vector.

The rest of this paper consists of five sections. Section~\ref{sec:M/G/1-type} provides preliminary results on M/G/1-type Markov chains. Section~\ref{sec:level increment_truncation} explains the LI truncation to compute the stationary distribution of an M/G/1-type Markov chain. Section~\ref{sec:convergence} shows the error caused by the LI truncation approximation converges, in total variation norm, to zero as the truncation parameter going to infinity. Section~\ref{sec:Subgeometric_convergence_formula} presents a subgeometric convergence formula for the level-wise difference between the original stationary distribution and its LI truncation approximation. Section~\ref{sec:concluding} contains concluding remarks.

\section{M/G/1-type Markov Chains}
\label{sec:M/G/1-type}
This section is divided into three subsections. Section~\ref{subsec:notation} introduces some basic notation and conventions used in the present and subsequent sections. Section~\ref{subsec-definition} provides the definition of the M/G/1-type Markov chain together with a well-known sufficient condition for the existence of the unique stationary distribution.  Finally, Section~\ref{subsec:Ramaswami} describes Ramaswami's recursion for the stationary distribution in M/G/1-type Markov chains.

\subsection{Notation and conventions}
\label{subsec:notation}

We define some sets of numbers. Let $\bbZ$ denote the set of all integers, and let
\begin{alignat*}{2}
\bbZ_{[k,\ell]} &= \{n \in \bbZ: k \le n \le \ell\},
&\quad
\bbZ_{\geqslant k} &= \{n \in \bbZ: n \ge k\},
\quad k,\ell\in\bbZ,
\\
\bbZ_+ &= \bbZ_{\geqslant 0}=\{n \in \bbZ: n \ge 0\},
&\quad
\bbN &= \bbZ_{\geqslant 1}=\{n \in \bbZ: n \ge 1\},
\\
\mathbb{M}_0 &= \{1,2,\ldots, M_0\},
&\qquad
\mathbb{M}_1 &= \{1,2,\ldots,M_1\},
\end{alignat*}
where $M_0,M_1 \in \bbN$.

We then introduce basic definitions and conventions on vectors and matrices.
Let $\vc{e}$ and $\vc{0}$ denote the column vectors of ones and zeros, respectively. Let $\vc{I}$ and $\vc{O}$ denote the identity matrix and the zero matrix, respectively. These vectors and matrices have appropriate dimensions (sizes) depending on where they are used. Furthermore, let $(\,\cdot\,)_{i,j}$ (resp.\ $(\,\cdot\,)_i$) denote the $(i,j)$-th (resp.\ $i$-th) element of the matrix (resp.\ vector) in the parentheses. For any matrix $\vc{X} \ge \vc {O}$ (resp.\ vector $\vc{x} \ge \vc{0}$), the notation $\vc{X} < \infty$ (resp.\ $\vc{x} < \infty$) means that every element of $\vc{X}$ (resp.\ $\vc{x}$) is finite. Finally, for any matrix function $\vc{Z}(\,\cdot\,)$ and scalar function $f(\,\cdot\,)$ on $(-\infty,\infty)$, the notation $\vc{Z}(x) = \ol{\bcal{O}}(f(x))$ means that
\begin{align*}
\limsup_{x\to\infty} { \sup_i \sum_{j}|(\vc{Z}(x))_{i,j}| \over f(x)}
&< \infty.
\end{align*}

\subsection{Definition and basic assumptions of the M/G/1-type Markov chain}\label{subsec-definition}

We define the M/G/1-type Markov chain. Let $\{(X_n, J_n); n \in \mathbb{Z}_+\}$ denote a discrete-time Markov chain on
state space $\bbS := \bigcup_{k=0}^\infty\{k\} \times \mathbb{M}_{k \vmin 1}$, where $x \vmin y = \min(x,y)$ for $x,y \in (-\infty,\infty)$. Let $\vc{P}$ denote the transition probability matrix
of the Markov chain $\{(X_n, J_n)\}$, and assume that $\vc{P}$ is a stochastic matrix with an upper block-Hessenberg structure:
\begin{equation}
\vc{P}=
\bordermatrix{
& \bbL_0
& \bbL_1
& \bbL_2
& \bbL_3
& \cdots \cr
\bbL_0
& \vc{B}(0)
& \vc{B}(1)
& \vc{B}(2)
& \vc{B}(3)
& \cdots \cr
\bbL_1
& \vc{B}(-1)
& \vc{A}(0)
& \vc{A}(1)
& \vc{A}(2)
& \cdots \cr
\bbL_2
& \vc{O}
& \vc{A}(-1)
& \vc{A}(0)
& \vc{A}(1)
& \cdots \cr
\bbL_3
& \vc{O}
& \vc{O}
& \vc{A}(-1)
& \vc{A}(0)
& \cdots
\cr
~\vdots
& \vdots
& \vdots
& \vdots
& \vdots
& \ddots
},
\label{defn-P}
\end{equation}
where $\bbL_k:= k \times \mathbb{M}_{k \vmin 1}$ is called level $k$ and an element $j \in \bbL_k$ is called {\it phase} $j$ of level $k$. The Markov chain $\{(X_n, J_n)\}$ is referred to as an M/G/1-type Markov chain \cite{Neut89}.

The following is the fundamental assumption of this paper, which is assumed throughout the paper unless otherwise stated.
\begin{ASSU} \label{assum1}
(i) The stochastic matrix $\vc{P}$ is irreducible;
(ii) $\vc{A}:=\sum_{k=-1}^\infty \vc{A}(k)$ is an irreducible stochastic matrix;
(iii) $\ol{\vc{m}}_B:=\sum_{k = 1}^\infty k\vc{B}(k)\vc{e} < \infty$; and
(iv) $\sigma := \vc{\varpi} \ol{\vc{m}}_A < 0$, where $\ol{\vc{m}}_A = \sum_{k = -1}^\infty k\vc{A}(k)\vc{e}$ and $\vc{\varpi}$ denotes the unique stationary distribution vector of $\vc{A}$.
\end{ASSU}

Assumption~\ref{assum1} ensures that $\vc{P}$ is irreducible and positive recurrent and thus it has the unique stationary distribution vector, denoted by $\vc{\pi} = (\pi(k, i))_{(k,i) \in \bbS}$ \cite[Chapter~XI, Proposition~3.1]{Asmu03}.
The subvectors $\{\vc{\pi}(k) := (\pi(k,i))_{i \in \bbM_{k\vmin 1}}$; $k \in \bbZ_+\}$ satisfy Ramaswami's recursion \cite{Rama88}, which is described in Section~\ref{subsec:Ramaswami}.

Assumption~\ref{assum1} also enables us to define the sequences of matrices that play an important role in deriving the main results in Section~\ref{sec:Subgeometric_convergence_formula}. Let $\{\ool{\vc{A}}(k); k\in\bbZ_{\geqslant-3}\}$, and $\{\ool{\vc{B}}(k); k\in\bbZ_{\geqslant-1}\}$ denote sets of matrices such that
\begin{alignat*}{2}
\ool{\vc{A}}(k)
&= \sum_{\ell = k + 1}^\infty \ol{\vc{A}}(\ell),
& \qquad
\ool{\vc{B}}(k)
&= \sum_{\ell = k + 1}^\infty \ol{\vc{B}}(\ell),
\end{alignat*}
where
\begin{alignat}{2}
\ol{\vc{A}}(k)
&= \sum_{\ell=k+1}^{\infty} \vc{A}(\ell),
& \qquad
\ol{\vc{B}}(k)
&= \sum_{\ell=k+1}^{\infty} \vc{B}(\ell).
\label{defn-ol{A}(k)}
\end{alignat}
Note that $\sum_{k=1}^{\infty}k\vc{A}(k) < \infty$ and $\sum_{k=1}^{\infty}k\vc{B}(k) < \infty$ if and only if $\sum_{k=0}^{\infty}\ol{\vc{A}}(k) < \infty$ and $\sum_{k=0}^{\infty}\ol{\vc{B}}(k) < \infty$, respectively. Therefore, $\ol{\vc{A}}(k) < \infty$ for $k \in \bbZ_{\geqslant -2}$ and $\ol{\vc{B}}(k) < \infty$ for $k \in \bbZ_+$ under Assumption~\ref{assum1} (especially the conditions (iii) and (iv)).

\subsection{Ramaswami's recursion}
\label{subsec:Ramaswami}

In this subsection, we describe Ramaswami's recursion for the stationary distribution in M/G/1-type Markov chains. To this end, we first introduce the $G$-matrix, the transition probability matrix over the first passage time to one level below in non-boundary levels. Using the $G$-matrix, we define probabilistically interpretable matrices, and then describe Ramaswami's recursion.

We begin with defining the $G$-matrix. Let $\vc{G}:=(G_{i, j})_{(i,j)\in(\bbM_1)^2}$ denote an $M_1 \times M_1$ matrices such that
\[
G_{i, j} = \PP\left((X_{T_1}, J_{T_1})
= (1, j) \mid (X_0, J_0) = (2, i)\right), \qquad
i, j \in \mathbb{M}_1,
\]
where $T_k= \inf \{n \in \mathbb{N}: (X_n,J_n) \in \bbL_k\}$ for $k \in \bbZ_+$. The matrix $\vc{G}$ is called the {\it $G$-matrix}, and is stochastic due to the conditions (ii) and (iv) of Assumption \ref{assum1} \cite[Theorem 2.3.1]{Neut89}. The stochastic matrix $\vc{G}$ has exact only one communication class \cite[Proposition~2.1]{Kimu10} and thus the unique stationary probability vector, denoted by $\vc{g}$.

The matrix $\vc{G}$ is generated by a recursion with infinite (computational) complexity. As is well known, $\vc{G}$ is the minimum nonnegative solution of the following matrix equation (see, e.g., \cite[Eq.~(2.3.3) and Theorem 2.2.2]{Neut89}):
\begin{equation*}
\mathscr{G} = \sum_{m = 0}^\infty \vc{A}(m - 1)\mathscr{G}^m.
\end{equation*}
Thus, we have a sequence $\{\vc{G}_n; n\in\bbZ_+\}$ converging to $\vc{G}$ from below (see, e.g., \cite{Gun89}):
\begin{subnumcases}{\label{recursion-G}}
{}
\vc{G}_0 = \vc{O},
\\
\vc{G}_n = \sum_{m = 0}^\infty \vc{A}(m - 1) (\vc{G}_{n-1})^m,
& \mbox{$n \in \bbN$}.
\end{subnumcases}

To describe Ramaswami's recursion, we introduce two matrices $\vc{\varPhi}(0)$ and $\vc{K}$. We first define $\vc{\varPhi}(0)$ as
\begin{alignat*}{2}
\vc{\varPhi}(0)
&= \sum _{m = 0}^\infty \vc{A}(m)\vc{G}^m.
\end{alignat*}
The definition of $\vc{G}$ leads to the interpretation:
\begin{alignat*}{2}
(\vc{\varPhi}(0))_{i, j}
&= \PP(J_{T_k} = j, T_k < T_{k-1} \mid X_0 = k, J_0 = i),
&\qquad i,j &\in \bbM_1.
\end{alignat*}
The irreducibility of $\vc{P}$ implies that the Markov chain $\{(X_n,J_n)\}$ does not keep staying in one level and thus $\sum_{n=0}^{\infty}(\vc{\varPhi}(0))^n=(\vc{I}-\vc{\varPhi}(0))^{-1}$. We then define $\vc{K}$ as
\begin{alignat*}{2}
\vc{K} &= \vc{B}(0) + \sum_{m = 1}^\infty \vc{B}(m)\vc{G}^{m-1}
(\vc{I} - \vc{\varPhi}(0))^{-1}\vc{B}(-1).
\end{alignat*}
By definition,
\begin{alignat*}{2}
(\vc{K})_{i, j}
&= \PP(J_{T_0} = j,  \mid X_0 = 0, J_0 = i),
&\qquad i,j &\in \bbM_0.
\end{alignat*}
Since $\{(X_n,J_n)\}$ is irreducible and positive recurrent (due to Assumption~\ref{assum1}), $\vc{K}$ is irreducible and stochastic, and thus it has the unique stationary probability, denoted by $\vc{\kappa}$.

We are ready to present Ramaswami's recursion for the stationary distribution vector $\vc{\pi}=(\vc{\pi}(0),\vc{\pi}(1),\dots)$.
\begin{PROP}[{}{\cite{Rama88},~\cite[Theorem 3.1]{Sche90}}]\label{prop-Rama88}
\begin{align*}
\vc{\pi}(0)
&= \vc{\kappa}
\Biggl[
1 + {\vc{\kappa} \over - \sigma}
\Biggl\{
\ol{\vc{m}}_B
+ \left(
\sum_{m = 1}^\infty \vc{B}(m)(\vc{I} - \vc{G}^m)
\right)(\vc{I} - \vc{A} + \vc{e}\vc{\varpi})^{-1}\ol{\vc{m}}_A
\Biggr\}
\Biggr]^{-1},
\\
\vc{\pi}(k)
&= \vc{\pi}(0)\vc{R}_0(k)
+ \sum_{\ell=1}^{k-1}\vc{\pi}(\ell)\vc{R}(k-\ell),
\qquad
k\in\bbN,
\end{align*}
where
\begin{alignat}{2}
\vc{R}(k)
&=
\dm\sum_{m=0}^{\infty}
\vc{A}(k+m)\vc{G}^m(\vc{I}-\vc{\varPhi}(0))^{-1}, & \qquad k &\in \bbN,
\label{defn-R(k)}
\\
\vc{R}_0(k)
&=
\sum_{m=0}^{\infty}
\vc{B}(k+m)\vc{G}^m(\vc{I}-\vc{\varPhi}(0))^{-1}, & \qquad k &\in \bbN.
\label{defn-R_0(k)}
\end{alignat}
\end{PROP}

\section{The LI Truncation Approximation}
\label{sec:level increment_truncation}

This section introduces the level-increment (LI) truncation approximation to the M/G/1-type Markov chain. The LI truncation approximation is often used to implement Ramaswami's recursion presented in Proposition~\ref{prop-Rama88}.
In the following, we explain the motivation of considering the LI truncation approximation and then define it.

Implementing Ramaswami's recursion requires the truncation of the infinite sequences $\{\vc{A}(k); k \in \bbZ_{\geqslant -1}\}$ and $\{\vc{B}(k); k \in \bbZ_+\}$. According to (\ref{defn-R(k)}) and (\ref{defn-R_0(k)}),  the main components $\vc{R}(k)$ and $\vc{R}_0(k)$ of Ramaswami's recursion are  computed through the infinite sums involving the infinite sequences $\{\vc{A}(k)\}$ and $\{\vc{B}(k)\}$. The infinite sum of $\{\vc{A}(k)\}$ is also included by the recursion (\ref{recursion-G}) of $\vc{G}$.

A standard way to run Ramaswami's recursion is to {\it truncate}
the infinite sequences $\{\vc{A}(k)\}$ and $\{\vc{B}(k)\}$ at some $k=N$, that is, to replace the two sequences with substantially finite sequences $\{\vc{A}^{(N)}(k)\}$ and $\{\vc{B}^{(N)}(k)\}$:
\begin{eqnarray}
\vc{A}^{(N)}(k)
&=&
\left\{
\begin{array}{ll}
\vc{A}(k),			& k \in \bbZ_{[-1,N-1]},
\\
\ol{\vc{A}}(N-1), 	& k=N,
\\
\vc{O}, 			& k \in \bbZ_{\geqslant N+1},
\end{array}
\right.
\label{defn-A^(N)(k)}
\\
\vc{B}^{(N)}(k)
&=&
\left\{
\begin{array}{ll}
\vc{B}(k),			& k \in \bbZ_{[-1,N-1]},
\\
\ol{\vc{B}}(N-1), 	& k=N,
\\
\vc{O}, 			& k \in \bbZ_{\geqslant N+1}.
\end{array}
\right.
\label{defn-B^(N)(k)}
\end{eqnarray}

We refer to the above-mentioned truncation as the {\it level-increment (LI) truncation}. This is because replacing $\vc{A}(k)$ and $\vc{B}(k)$, $k\in \bbZ_{\geqslant N+1}$ with zero matrices can be interpreted as truncating the level increment per transition at the upper limit $N$. Indeed, the LI truncation transforms the transition probability matrix $\vc{P}$ into another M/G/1-type stochastic matrix
\begin{align}
\vc{P}^{(N)}
:=
\bordermatrix{
& \bbL_0
& \bbL_1
& \bbL_2
& \bbL_3
& \cdots \cr
\bbL_0
& \vc{B}^{(N)}(0)
& \vc{B}^{(N)}(1)
& \vc{B}^{(N)}(2)
& \vc{B}^{(N)}(3)
& \cdots \cr
\bbL_1
& \vc{B}^{(N)}(-1)
& \vc{A}^{(N)}(0)
& \vc{A}^{(N)}(1)
& \vc{A}^{(N)}(2)
& \cdots \cr
\bbL_2
& \vc{O}
& \vc{A}^{(N)}(-1)
& \vc{A}^{(N)}(0)
& \vc{A}^{(N)}(1)
& \cdots \cr
\bbL_3
& \vc{O}
& \vc{O}
& \vc{A}^{(N)}(-1)
& \vc{A}^{(N)}(0)
& \cdots
\cr
~\vdots
& \vdots
& \vdots
& \vdots
& \vdots
& \ddots
},
\label{defn-P^{(N)}}
\end{align}
and a Markov chain driven by $\vc{P}^{(N)}$ has level jumps of at most $N$.

We note that performing Ramaswami's recursion by the LI truncation is equivalent to compute a stationary distribution of the M/G/1-type stochastic matrix $\vc{P}^{(N)}$. The following proposition ensures that $\vc{P}^{(N)}$ has the unique stationary distribution under Assumption~\ref{assum1}.
\begin{PROP}\label{lem-communication-P^{(N)}}
If Assumption~\ref{assum1} holds, then, for any $N \in \bbN$, $\vc{P}^{(N)}$ is positive Harris recurrent having a Harris recurrent set $\bbL_0$ (see, e.g., \cite[Chapter 9]{Meyn09}).
\end{PROP}
\proof
We consider a Markov chain $\{(X_n^{(N)},J_n^{(N)}); n\in\bbZ_+\}$ with transition probability matrix $\vc{P}^{(N)}$ on the same probability space as the M/G/1-type Markov chain $\{(X_n,J_n); n\in\bbZ_+\}$. According to the definition (\ref{defn-P^{(N)}}) of $\vc{P}^{(N)}$, we assume without loss of generality that
\begin{subequations}\label{eqn:210302-01}
\begin{align}
X_n^{(N)} &=
\left\{
\begin{array}{ll}
X_0, & n=0,
\\
X_{n-1}^{(N)} + \min(X_n - X_{n-1},N), &
n=1,2,\dots,T_0^{(N)},
\end{array}
\right.
\\
J_n^{(N)} &= J_n,
\quad n =0,1,\dots,T_0^{(N)}-1,
\end{align}
\end{subequations}
where $T_0^{(N)} = \inf\{n\in\bbN: X_n^{(N)} = 0\}$. Equation (\ref{eqn:210302-01}) implies that $\{(X_n^{(N)},J_n^{(N)})\}$ can reach non-boundary levels $\bbS\setminus \bbL_0$ from any state in $\bbL_0$. Equation (\ref{eqn:210302-01}) also implies that
\begin{subequations}\label{eqn:210114-01}
\begin{alignat}{2}
X_n^{(N)} &\le X_n,&\quad n &= 0,1,\dots,T_0^{(N)},
\\
J_n^{(N)} &= J_n,  &\quad n &= 0,1,\dots,T_0^{(N)}-1.
\end{alignat}
\end{subequations}
Equation (\ref{eqn:210114-01}) shows that $\{(X_n,J_n)\}$ and thus $\{(X_n^{(N)},J_n^{(N)})\}$ can reach $\bbL_0$ from any state in non-boundary levels $\bbS\setminus \bbL_0$. Note here that $\{(X_n,J_n)\}$ is irreducible and positive recurrent under Assumption~\ref{assum1}. Therefore, from (\ref{eqn:210302-01})  and (\ref{eqn:210114-01}), we have
\begin{align*}
\EE\left[T_0^{(N)}
\mid (X_0^{(N)},J_0^{(N)}) = (k,i) \right]
&\le
\EE\left[T_0
\mid (X_0,J_0) = (k,i) \right] < \infty
\quad \mbox{for all $(k,i) \in \bbS$}.
\end{align*}
The proof is completed.
\qed

Based on Proposition~\ref{lem-communication-P^{(N)}}, we denote by $\vc{\pi}^{(N)}$ the unique stationary distribution of $\vc{P}^{(N)}$, and then refer to $\vc{\pi}^{(N)}$ and $\vc{P}^{(N)}$ as the level-increment (LI) truncation approximations to $\vc{\pi}$ and $\vc{P}$, respectively.

\section{Convergence of the LI Truncation Approximation}
\label{sec:convergence}

In this section, we first derive a difference formula for the original stationary distribution $\vc{\pi}$ and its LI truncation approximation $\vc{\pi}^{(N)}$. With the difference formula, we show the basic convergence of $\vc{\pi}^{(N)}$ to $\vc{\pi}$ as $N \to \infty$.

We begin with the definition of a matrix associated with the deviation matrix. Fix $(k_*, i_*)\in \bbS$ arbitrarily, and let $\vc{H} := (H(k,i ; \ell,j))_{(k,i ; \ell,j)\in\bbS^2}$ denote a matrix such that
\begin{equation*}
H(k,i;\ell,j) = \EE_{(k,i)}\!
\left[
\sum_{\nu=0}^{T_{(k_*,i_*)}-1} \dd{1}((X_{\nu},J_{\nu})=(\ell,j))
\right]
- \pi(\ell,j)\EE_{(k,i)}
\left[ T_{(k_*,i_*)} \right],
\end{equation*}
where $T_{(k_*,i_*)} = \inf\{\nu \in \bbN: (X_{\nu},J_{\nu}) = (k_*,i_*)\}$,
and where
\[
\EE_{(k,i)}[\,\cdot\,]=\EE[\,\cdot\mid(X_0,J_0)=(k,i)],\quad(k,i) \in \bbS.
\]
The matrix $\vc{H}$ is well-defined under Assumption \ref{assum1}.

\begin{REM}
The matrix $\vc{H}$ depends on the state $(k_*, i_*)$, but essentially the choice of $(k_*, i_*)$ does not affect the subsequent analysis, and moreover the subgeometric convergence formula presented in Theorem~\ref{th:subexp} (the main result of this paper) does not depend on $\vc{H}$. Thus, we do not parameterize $\vc{H}$ with $(k_*, i_*)$ for simplicity.
\end{REM}

Using $\vc{H}$, we present a difference formula for $\vc{\pi}^{(N)}$ and $\vc{\pi}$. Note that $\vc{H}$ is equivalent to $\wt{\vc{H}}$ defined in the proof of \cite[Theorem~4.1]{Masu21-M/G/1-Subexp} and $\vc{H}$ satisfies the following:
\[
(\vc{I} - \vc{P})\vc{H} = \vc{I} - \vc{e}\vc{\pi}.
\]
Therefore, we have a difference formula of the same form as \cite[Lemma~5.1]{Masu21-M/G/1-Subexp}:
\begin{equation*}
\vc{\pi}^{(N)}-\vc{\pi}=\vc{\pi}^{(N)}(\vc{P}^{(N)}-\vc{P})\vc{H},
\end{equation*}
which leads to
\begin{equation}
\vc{\pi}^{(N)}(k) - \vc{\pi}(k)
= \sum_{\ell=0}^{\infty} \vc{\pi}^{(N)}(\ell)
\sum_{n=0}^{\infty} \vc{\varDelta}^{(N)}(\ell;n) \vc{H}(n;k),
\qquad k \in \mathbb{Z}_+,
\label{eq:piN-pi1-add}
\end{equation}
where $\vc{H}(k; \ell)$ and $\vc{\varDelta}^{(N)}(k; \ell)$, $k,\ell\in\bbZ_+$, denote the $(k, \ell)$-blocks of $\vc{H}$ and $\vc{\varDelta}^{(N)}:=\vc{P}^{(N)} - \vc{P}$, respectively.

To present a more detailed expression for $\vc{\pi}^{(N)}(k) - \vc{\pi}(k)$ than (\ref{eq:piN-pi1-add}), we introduce a probabilistically interpretable matrix $\vc{F}_+$: Let $\vc{F}_+ = (F_+(m,i; k,j))_{(m,i; k,j) \in (\bbS)^2}$ denote a nonnegative matrix such that
\begin{equation*}
F_+(m,i;k,j)
=\EE_{(m,i)}\!
\left[
\sum_{\nu=0}^{T_0-1} \dd{1}((X_{\nu},J_{\nu})=(k,j))
\right],
\end{equation*}
where $T_0 =\inf_{j\in\bbM_0} T_{(0,j)} = \inf\{n\in\bbN: X_n=0\}$.

The following lemma presents a detailed expression of the level-wise difference  $\vc{\pi}^{(N)}(k) - \vc{\pi}(k)$.
\begin{LEM}\label{lem:difference_formula}
If Assumption \ref{assum1} holds, then, for any $N \in \bbN$,
\begin{align}
&\vc{\pi}^{(N)}(k) - \vc{\pi}(k)
\nonumber
\\
&~~= \vc{\pi}^{(N)}(0)
\left[
{1 \over -\sigma}
\ool{\vc{B}}(N - 1)\vc{e}
\vc{\pi}(k)
+
\sum_{n=N+1}^{\infty}
\vc{B}(n) ( \vc{G}^{N-k} - \vc{G}^{n-k} )
\vc{F}_+(k;k)
\right.
\nonumber
\\
&\qquad\qquad\qquad\qquad\quad
\Biggl.
+
\sum_{n=N+1}^{\infty}
\vc{B}(n)( \vc{G}^{N-1} - \vc{G}^{n-1} )
\vc{S}(k)
\Biggr]
\nonumber
\\
&~~~~
+ \sum_{\ell=1}^{\infty}
\vc{\pi}^{(N)}(\ell)
\left[
{1 \over -\sigma}
\ool{\vc{A}}(N - 1)\vc{e}
\vc{\pi}(k)
+ \sum_{n=N+1}^{\infty}
\vc{A}(n) (\vc{G}^{N + \ell - k} - \vc{G}^{n + \ell - k})
\vc{F}_+(k;k)
\right.
\qquad
\nonumber
\\
&\qquad\qquad\qquad\qquad\quad
\Biggl.
+ \sum_{n=N+1}^{\infty}
\vc{A}(n) (\vc{G}^{N + \ell - 1} - \vc{G}^{n + \ell - 1})
\vc{S}(k)
\Biggr],\quad k \in \bbZ_{[0,N]},
\label{eq:piN(k)-pi(k)}
\end{align}
where
\begin{alignat}{2}
\vc{S}(k) &= (\vc{I}-\vc{\Phi}(0))^{-1}\vc{B}(-1)\vc{H}(0; k)
+ \vc{G}(\vc{I} -\vc{A} -\ol{\vc{m}}_A\vc{g})^{-1}\vc{e}\vc{\pi}(k),
&\quad
k &\in \mathbb{Z}_+.
\label{eq:defS}
\end{alignat}
\end{LEM}

\proof
See Appendix~\ref{proof:lem:difference_formula}.
\qed

Using Lemma \ref{lem:difference_formula}, we show that $\{\vc{\pi}^{(N)}; N \in \mathbb{N}\}$ converges to $\vc{\pi}$ in {\it total variation norm}. Total variation norm is defined as follows:  For any $\vc{x}$, $|\vc{x}|$ denotes a nonnegative vector (or matrix) by taking the absolute values of the elements of $\vc{x}$, and $\|\vc{x}\|$ denotes the total variation norm of the vector $\vc{x}$, i.e., $\|\vc{x}\| = |\vc{x}| \vc{e}$.

The following result is the goal of this section.
\begin{THM}\label{th:convergence_tvn}
If Assumption \ref{assum1} holds, then
\begin{equation}
\lim_{N \to \infty} \|\vc{\pi}^{(N)} - \vc{\pi}\| = 0.
\label{eq:th-convergence_tvn}
\end{equation}
\end{THM}
\proof
To prove (\ref{eq:th-convergence_tvn}), it suffices to show that
\begin{equation}
\lim_{N \to \infty} \vc{\pi}^{(N)}(k) = \vc{\pi}(k)
\quad \mbox{for all $k \in\bbZ_+$}.
\label{eq:th-convergence}
\end{equation}
Indeed, $\|\vc{\pi}^{(N)} - \vc{\pi}\| \le 2$ for all $N \in \bbN$, and therefore it follows from (\ref{eq:th-convergence}) and the dominated convergence that
\begin{align*}
\lim_{N\to\infty} \|\vc{\pi}^{(N)} - \vc{\pi}\|
&= \lim_{N\to\infty} \sum_{k=0}^\infty |\vc{\pi}^{(N)}(k) - \vc{\pi}(k)|\vc{e}
\\
&= \sum_{k=0}^\infty \lim_{N\to\infty} |\vc{\pi}^{(N)}(k) - \vc{\pi}(k)|\vc{e}
\\
&= 0,
\end{align*}
which yields \eqref{eq:th-convergence_tvn}.

In what follows, we prove (\ref{eq:th-convergence}). By definition,
\begin{subequations}\label{eqn:210211-01}
\begin{align}
\vc{\pi}^{(N)}(0) &\le \vc{e}^\top,
\quad
\sum_{\ell=1}^{\infty}\vc{\pi}^{(N)}(\ell)  \le \vc{e}^\top,
\\
|\vc{G}^m - \vc{G}^n| &\le 2\vc{e}\vc{e}^{\top}\quad \mbox{for all $m,n \in \bbZ_+$}.
\end{align}
\end{subequations}
Using these inequalities and (\ref{eq:piN(k)-pi(k)}), we obtain
\begin{align}
&|\vc{\pi}^{(N)}(k) - \vc{\pi}(k)|
\nonumber
\\
&~~\le \vc{\pi}^{(N)}(0)
\left[
{1 \over -\sigma}
\ool{\vc{B}}(N - 1)\vc{e}
\vc{\pi}(k)
+
\sum_{n=N+1}^{\infty}
\vc{B}(n) \left| \vc{G}^{N-k} - \vc{G}^{n-k} \right|
\vc{F}_+(k;k)
\right.
\nonumber
\\
&\qquad\qquad\qquad\qquad\quad
\Biggl.
+
\sum_{n=N+1}^{\infty}
\vc{B}(n)\left| \vc{G}^{N-1} - \vc{G}^{n-1} \right|
\left|\vc{S}(k)\right|
\Biggr]
\nonumber
\\
&~~~~
+ \sum_{\ell=1}^{\infty}
\vc{\pi}^{(N)}(\ell)
\left[
{1 \over -\sigma}
\ool{\vc{A}}(N - 1)\vc{e}
\vc{\pi}(k)
+ \sum_{n=N+1}^{\infty}
\vc{A}(n)
\left|\vc{G}^{N + \ell - k} - \vc{G}^{n + \ell - k}\right|
\vc{F}_+(k;k)
\right.
\qquad
\nonumber
\\
&\qquad\qquad\qquad\qquad\quad
\Biggl.
+ \sum_{n=N+1}^{\infty}
\vc{A}(n) \left| \vc{G}^{N + \ell - 1} - \vc{G}^{n + \ell - 1}\right|
\left|\vc{S}(k)\right|
\Biggr]
\nonumber
\\
&~~\le \vc{e}^\top
\left[
{1 \over -\sigma}
\ool{\vc{B}}(N - 1)\vc{e}
\vc{\pi}(k)
+
2\ol{\vc{B}}(N)\vc{e}\vc{e}^\top
\left\{\vc{F}_+(k;k) + \left|\vc{S}(k)\right|
\right\}
\right]
\nonumber
\\
&~~~~
+ \vc{e}^\top
\left[
{1 \over -\sigma}
\ool{\vc{A}}(N - 1)\vc{e}
\vc{\pi}(k)
+ 2\ol{\vc{A}}(N)\vc{e}\vc{e}^\top
\left\{
\vc{F}_+(k;k) + \left|\vc{S}(k)\right|
\right\}
\right],
~~ k\in\bbZ_{[0,N]}.
\label{eq:piN(k)-pi(k)_2}
\end{align}
From Assumption~\ref{assum1}, we also have
\begin{alignat*}{2}
\lim_{N\to\infty}\ol{\vc{A}}(N)\vc{e} &= \vc{0}, &\quad
\lim_{N\to\infty}\ool{\vc{A}}(N)\vc{e} &= \vc{0},
\\
\lim_{N\to\infty}\ol{\vc{B}}(N)\vc{e} &= \vc{0}, &\quad
\lim_{N\to\infty}\ool{\vc{B}}(N)\vc{e} &= \vc{0}.
\end{alignat*}
Applying these equations to (\ref{eq:piN(k)-pi(k)_2}) yields
\[
\lim_{N\to\infty}|\vc{\pi}^{(N)}(k) - \vc{\pi}(k)| = \vc{0},
\qquad k \in \bbZ_+,
\]
and thus \eqref{eq:th-convergence} holds. The proof is completed.
\qed

\section{Subgeometric Convergence Formulas for the LI Truncation Approximation}
\label{sec:Subgeometric_convergence_formula}

This section presents subgeometric convergence formulas for $\vc{\pi}^{(N)}(k)-\vc{\pi}(k)$ under the following assumption.
\begin{ASSU}\label{assum3}
There exists a distribution function $F$ on $\bbZ_+$ such that $F$ is {\it long-tailed} (i.e., $F \in \calL$; see Definition~\ref{defn-long-tailed}~(i)) and
\[
\lim_{N \to \infty}
{
\ool{\vc{A}}(N)\vc{e}
\over
\ol{F}(N)
} = \vc{c}_A,
\quad
\lim_{N \to \infty}
{
\ool{\vc{B}}(N)\vc{e}
\over
\ol{F}(N)
} = \vc{c}_B,
\]
where $\vc{c}_A \ge \vc{0}$ and $\vc{c}_B \ge \vc{0}$ are $M_1$- and $M_0$-dimensional finite column vectors, respectively, and either of them is a non-zero vector.
\end{ASSU}

The following theorem is the main result of this paper.
\begin{THM}\label{th:subexp}
If Assumptions \ref{assum1} and \ref{assum3} hold,
then
\begin{subequations}\label{eqns:thm}
\begin{alignat}{2}
\lim_{N \to \infty}
{
\vc{\pi}^{(N)}(k) - \vc{\pi}(k)
\over
\ol{F}(N)
}
&=
{\vc{\pi}(0)\vc{c}_B + \ol{\vc{\pi}}(0)\vc{c}_A
\over
-\sigma
}\vc{\pi}(k)>\vc{0}, &\qquad k &\in \bbZ_+,
\label{eq:subexp-01}
\end{alignat}
and thus
\begin{alignat}{2}
\lim_{N \to \infty}
{1 \over \ol{F}(N)}
\left\|
{
\vc{\pi}^{(N)}(k) - \vc{\pi}(k)
\over
\vc{\pi}(k)\vc{e}
}
\right\|
&=
{\vc{\pi}(0)\vc{c}_B + \ol{\vc{\pi}}(0)\vc{c}_A
\over
-\sigma
}>0, &\qquad k &\in \bbZ_+,
\label{eq:subexp-norm-01}
\end{alignat}
where $\ol{\vc{\pi}}(k) = \sum_{\ell = k + 1}^\infty \vc{\pi}(\ell)$ for $k \in \mathbb{Z}_+$.
\end{subequations}

\end{THM}

\proof
See Appendix~\ref{proof:th:subexp}.
\qed

Theorem~\ref{th:subexp} shows the subgeometric convergence of the level-wise difference $\vc{\pi}^{(N)}(k) - \vc{\pi}(k)$ appears under Assumptions \ref{assum1} and \ref{assum3}, and this subgeometric convergence is connected with the level increment of the original chain $\{(X_n,J_n)\}$. To see this, we introduce a certain distribution associated with it, as in \cite{Masu21-M/G/1-Subexp}. Let $D$ denote a probability distribution (function) such that
\begin{align}
D(k)
&= \sum_{(\ell,i) \in \bbS} \pi(\ell,i)
\PP( \max(X_1-X_0,0) \le k \mid (X_0,J_0)=(\ell,i))
\nonumber
\\
&= \sum_{n=0}^k
\left[ \vc{\pi}(0) \vc{B}(n)\vc{e}
+ \ol{\vc{\pi}}(0) \cdot
\left\{ \vc{A}(n)\vc{e} + \delta_{n,0} \vc{A}(-1)\vc{e} \right\}
\right].
\label{defn:D}
\end{align}
The distribution $D$ is referred to as the {\it Nonnegative Level-increment-in-Steady-state (NLS) distribution}. Furthermore, let $D_I$ denote the integrated-tail distribution (or the equilibrium distribution) of the NLS distribution $D$, that is,
\begin{alignat}{2}
D_I(k)
&= {\sum_{\ell=0}^k (1 - D(\ell)) \over \sum_{\ell=1}^{\infty}\ell D(\ell)}
=
\sum_{\ell=0}^k
{
\vc{\pi}(0) \ol{\vc{B}}(\ell)\vc{e} + \ol{\vc{\pi}}(0) \ol{\vc{A}}(\ell)\vc{e}
\over
\vc{\pi}(0) \ol{\vc{m}}_{B} + \ol{\vc{\pi}}(0) \ol{\vc{m}}_{A}^+
},
\qquad k  \in \bbZ_+,
\label{defn:D_I(k)}
\end{alignat}
where the second equality is due to (\ref{defn:D}) and $\ol{\vc{m}}_{A}^+ :=\sum_{k=1}^{\infty}k\vc{A}(k)\vc{e}=\ol{\vc{m}}_{A} + \vc{A}(-1)\vc{e}$. It follows from (\ref{defn:D_I(k)}) and Assumption~\ref{assum3} that $\ol{D}_I:=1 - D_I$ satisfies
\begin{align}
\lim_{k\to\infty}{\ol{D}_I(k) \over \ol{F}(k)}
= {
\vc{\pi}(0) \vc{c}_{B} + \ol{\vc{\pi}}(0) \vc{c}_{A}
\over
\vc{\pi}(0) \ol{\vc{m}}_{B} + \ol{\vc{\pi}}(0) \ol{\vc{m}}_{A}^+
} \in (0,\infty).
\label{asymp:F_pi^I}
\end{align}
Therefore, Theorem~\ref{th:subexp} leads to the following result.
\begin{COR}\label{cor:subexp-norm}
If all the conditions of Theorem~\ref{th:subexp} are satisfied, then
\begin{subequations}\label{eqns:coro1}
\begin{alignat}{2}
\lim_{N\to\infty}
{
\vc{\pi}^{(N)}(k) - \vc{\pi}(k)
\over
\ol{D}_I(N)
}
&=
{\vc{\pi}(0) \ol{\vc{m}}_{B} + \ol{\vc{\pi}}(0) \ol{\vc{m}}_{A}^+
 \over -\sigma}
\vc{\pi}(k)>\vc{0}, &\qquad k &\in \bbZ_+,
\label{eq:subexp-02}
\\
\lim_{N\to\infty}
{1 \over
\ol{D}_I(N)
}
{
\left\|
\vc{\pi}^{(N)}(k) - \vc{\pi}(k)
\over
\vc{\pi}(k)\vc{e}
\right\|
}
&=
{\vc{\pi}(0) \ol{\vc{m}}_{B} + \ol{\vc{\pi}}(0) \ol{\vc{m}}_{A}^+
 \over -\sigma}>0, &\qquad k &\in \bbZ_+.
\label{eq:subexp-norm-02}
\end{alignat}
\end{subequations}

\end{COR}

\proof
Equation~(\ref{eq:subexp-01}) shows that for each $k \in \bbZ_+$ there exists some $N_k \in \bbN$ such that $\vc{\pi}^{(N)}(k) - \vc{\pi}(k) > \vc{0}$ for all $N \ge N_k$. Therefore,
\begin{align*}
\left\|
{
\vc{\pi}^{(N)}(k) - \vc{\pi}(k)
\over
\vc{\pi}(k)\vc{e}
}
\right\|
=
{
( \vc{\pi}^{(N)}(k) - \vc{\pi}(k) ) \vc{e}
\over
\vc{\pi}(k)\vc{e}
},\qquad N \ge N_k.
\end{align*}
Applying (\ref{eq:subexp-01}) to the above equation leads to (\ref{eq:subexp-norm-01}). Furthermore, combining (\ref{eq:subexp-norm-01}) with (\ref{asymp:F_pi^I}) yields (\ref{eq:subexp-norm-02}). The proof is completed.
\qed

Theorem~\ref{th:subexp} and Corollary~\ref{cor:subexp-norm} present subgeometric convergence formulas for the level-wise difference $\vc{\pi}^{(N)}(k) - \vc{\pi}(k)$. Equations~(\ref{eq:subexp-01}) and (\ref{eq:subexp-02}) show that, as the truncation parameter $N$ goes to $\infty$, the level-wise difference $\vc{\pi}^{(N)}(k) - \vc{\pi}(k)$ converges to zero at the same speed as subgeometric functions $\ol{F}$ and $\ol{D}_I$ (see Remark~\ref{rem:H=Lambda}). Moreover, \eqref{eq:subexp-norm-01} and \eqref{eq:subexp-norm-02} show that the total variation norm of the relative difference $\{\vc{\pi}^{(N)}(k) - \vc{\pi}(k) \} / \vc{\pi}(k)\vc{e}$ decays asymptotically at a rate independent of the level $k$.

Remarkably, Theorem~\ref{th:subexp} and Corollary~\ref{cor:subexp-norm} do not require the subexponentiality of $F$. Asymptotic analysis involving heavy-tailed distributions often assumes that a reference distribution specifying the decay rate, such as $F$, is subexponential. Indeed, Masuyama et al.\ \cite{Masu21-M/G/1-Subexp} presented similar subgeometric convergence formulas for the LCBA truncation approximation, though the formulas require the subexponentiality of $F$ (and thus $D_I$) and the finiteness of the second moment of level increments; more specifically,
\begin{enumerate}
\item $F \in \calS$, or equivalently, $D_I \in \calS$; and
\item $\sum_{k=1}^{\infty}k^2 \vc{A}(k)\vc{e} < \infty$ and $\sum_{k=1}^{\infty}k^2 \vc{B}(k)\vc{e} < \infty$.
\end{enumerate}
In short, the subgeometric convergence formulas for the LCBA truncation approximation are of the same forms as our ones, but they need stronger conditions.

We note that if $F \in \calS$ (as in \cite{Masu21-M/G/1-Subexp}) then the subgeometric convergence of the level-wise difference $\vc{\pi}^{(N)}(k) - \vc{\pi}(k)$ is connected with the tail decay of the stationary distribution of the original chain $\{(X_n,J_n)\}$.
\begin{COR}\label{cor:subexp}
Suppose that Assumptions \ref{assum1} and \ref{assum3} are satisfied. If $F \in \calS$, then
\begin{subequations}
\begin{alignat}{2}
\lim_{N \to \infty}
{
\vc{\pi}^{(N)}(k) - \vc{\pi}(k)
\over
\ol{\vc{\pi}}(N)\vc{e}
} &= \vc{\pi}(k), &\qquad k &\in\bbZ_+,
\label{eq:coro-subexp}
\\
\lim_{N \to \infty}
{1 \over \ol{\vc{\pi}}(N)\vc{e}}
\left\|
{
\vc{\pi}^{(N)}(k) - \vc{\pi}(k)
\over
\vc{\pi}(k)\vc{e}
}
\right\|
&=
1,
&\qquad k &\in\bbZ_+.
\label{eq:coro-subexp-norm}
\end{alignat}
\end{subequations}
\end{COR}
\proof
It follows from \cite[Theorem 3.1]{Masu16-ANOR} that if Assumptions \ref{assum1} and \ref{assum3} hold with $F \in \calS$ then
\begin{align*}
\lim_{N \to \infty}
{
\ol{\vc{\pi}}(N)
\over
\ol{F}(N)
}
=
{
\vc{\pi}(0)\vc{c}_B + \ol{\vc{\pi}}(0)\vc{c}_A
\over
-\sigma
}\vc{\varpi}.
\end{align*}
Using this equation and Theorem \ref{th:subexp}, we obtain
\begin{align*}
\lim_{N \to \infty}
\frac{\vc{\pi}^{(N)}(k) - \vc{\pi}(k)}{\ol{\vc{\pi}}(N)\vc{e}}
&= \lim_{N \to \infty}
\frac{\vc{\pi}^{(N)}(k) - \vc{\pi}(k)}{\ol{F}(N)} \cdot
\frac{\ol{F}(N)}{\ol{\vc{\pi}}(N)\vc{e}}
\\
&= \lim_{N \to \infty}
\frac{\vc{\pi}^{(N)}(k) - \vc{\pi}(k)}{\ol{F}(N)} \cdot
\lim_{N \to \infty} \frac{\ol{F}(N)}{\ol{\vc{\pi}}(N)\vc{e}}
\\
&= \frac{\vc{\pi}(0)\vc{c}_B + \ol{\vc{\pi}}(0)\vc{c}_A}{-\sigma}\vc{\pi}(k)
\cdot
\frac{-\sigma}{\vc{\pi}(0)\vc{c}_B + \ol{\vc{\pi}}(0)\vc{c}_A}
\\
&= \vc{\pi}(k),\qquad k\in\bbZ_+,
\end{align*}
which shows that (\ref{eq:coro-subexp}) holds. Furthermore,  (\ref{eq:coro-subexp-norm}) follows from (\ref{eq:coro-subexp}).
\qed

\begin{REM}
Our formulas (\ref{eq:coro-subexp}) and (\ref{eq:coro-subexp-norm}) in
Corollary~\ref{cor:subexp} are in the same form as the corresponding ones presented in \cite[Corollaries~5.13 and 5.14]{Masu21-M/G/1-Subexp}, though our formulas do not necessarily require $\sum_{k=1}^{\infty}k^2 \vc{A}(k)\vc{e} < \infty$ and $\sum_{k=1}^{\infty}k^2 \vc{B}(k)\vc{e} < \infty$. Therefore,
our formulas (\ref{eq:coro-subexp}) and (\ref{eq:coro-subexp-norm}) holds when  $F$ is a Pareto distribution even if its shape parameter is not more than two (see Example~\ref{exa:Pareto}).
\end{REM}

\section{Concluding Remarks}\label{sec:concluding}

We have derived several types of level-wise subgeometric convergence formulas for the level-increment (LI) truncation approximation $\vc{\pi}^{(N)}$ to the stationary distribution $\vc{\pi}$ of an M/G/1-type Markov chain. Although our subgeometric convergence formulas for the LI truncation approximation are the same as the corresponding ones for the last-column-block-augmented (LCBA) truncation approximation in \cite{Masu21-M/G/1-Subexp}, the former ones hold under weaker conditions than the latter ones.

It should be noted that, just because the formulas for the LI truncation approximation are the same type as the ones for the LCBA truncation approximation, does not necessarily mean that the accuracy of both is about the same. Provided that their truncation parameters are set to be the same value $N$, the LI truncation of the M/G/1-type stochastic matrix $\vc{P}$ is closer to the original stochastic matrix than the LCBA truncation of $\vc{P}$. Thus, it is likely that the LI truncation approximation to the original stationary distribution vector $\vc{\pi}$ is more accurate than the LCBA truncation approximation to $\vc{\pi}$. However, the difference between the two approximations is expected to be relatively negligible compared to the decay rate of the tail distributions $\ol{F}$ and thus $\ol{D}_I$. Evaluating this net difference between the two approximations is an interesting future task.

There are two other interesting problems associated with this study. One is to derive a subgeometric convergence formula for the total variation of the whole (not level-wise) difference between the stationary distribution and its LI truncation approximation. If such a uniform convergence was shown, then that result could help us determine the truncation parameter $N$ to meet a given error tolerance for computing the whole stationary distribution. The other is to derive geometric convergence formulas for the whole and/or level-wise difference between the stationary distribution and its LI truncation approximation. If this second problem were solved, we could see in what case it is easy to compute approximately the stationary distribution $\vc{\pi}$ of an M/G/1-type Markov chain.

\appendix

\section{Proof of Lemma~\ref{lem:difference_formula}}
\label{proof:lem:difference_formula}
First of all, we describe $\vc{\pi}^{(N)}(k) - \vc{\pi}(k)$ by the block component matrices $\vc{A}(k)$ and $\vc{B}(k)$ of $\vc{P}$ together with those of $\vc{H}$. From \eqref{defn-P} and \eqref{defn-A^(N)(k)}--\eqref{defn-P^{(N)}}, we have
\begin{equation}
\vc{\varDelta}^{(N)}(\ell;n)
=
\left\{
\begin{array}{lll}
\ol{\vc{B}}(N), & \ell = 0, & n = N,
\\
-\vc{B}(n), & \ell = 0, & n=N+1,N+2,\dots,
\\
\ol{\vc{A}}(N), & \ell \in \bbN, & n = N + \ell,
\\
-\vc{A}(n-\ell),  & \ell \in \bbN, & n=N + \ell + 1,N + \ell + 2,\dots,
\\
\vc{O}, & \mathrm{otherwise}.
\end{array}
\right.
\label{eq:PN-P}
\end{equation}
Substituting \eqref{eq:PN-P} into \eqref{eq:piN-pi1-add} yields
\begin{align}
&\vc{\pi}^{(N)}(k) - \vc{\pi}(k)
\nonumber
\\
&\quad = \vc{\pi}^{(N)}(0)\left\{\ol{\vc{B}}(N)\vc{H}(N; k)
- \sum_{n = N + 1}^{\infty} \vc{B}(n)\vc{H}(n; k)\right\}
\nonumber
\\
&\quad\quad~~~ + \sum_{\ell = 1}^{\infty} \vc{\pi}^{(N)}(\ell)
\left\{\ol{\vc{A}}(N)\vc{H}(N + \ell; k)
- \sum_{n=N+1}^{\infty} \vc{A}(n)\vc{H}(n + \ell; k)\right\}
\nonumber
\\
&\quad = \vc{\pi}^{(N)}(0)
\sum_{n = N + 1}^{\infty} \vc{B}(n)
\left\{ \vc{H}(N; k) - \vc{H}(n; k) \right\}
\nonumber
\\
&\quad\quad~~~ + \sum_{\ell = 1}^{\infty} \vc{\pi}^{(N)}(\ell)
\sum_{n = N + 1}^{\infty} \vc{A}(n)
\left\{\vc{H}(N + \ell; k) - \vc{H}(n + \ell; k)\right\},
\quad k\in\bbZ_{[0,N]}.
\label{eq:piN-pi2}
\end{align}

To derive (\ref{eq:piN(k)-pi(k)}) from (\ref{eq:piN-pi2}), we rearrange the two terms:
\begin{align*}
\sum_{n = N + 1}^{\infty} \vc{B}(n)
\left\{ \vc{H}(N; k) - \vc{H}(n; k) \right\}
\quad
\mbox{and}
\quad
\sum_{n = N + 1}^{\infty} \vc{A}(n) \left\{\vc{H}(N + \ell; k) - \vc{H}(n + \ell; k)\right\}.
\end{align*}
It follows from \cite[Theorem 9]{Zhao03} and \cite[Remark 4.8]{Masu21-M/G/1-Subexp} that, for $m > k$,
\begin{align}
\vc{H}(m;k)
&=\vc{F}_+(m; k)
+ \vc{G}^{m-1}(\vc{I}-\vc{\Phi}(0))^{-1}\vc{B}(-1)\vc{H}(0;k)
- \vc{u}(m)\vc{\pi}(k),
\nonumber
\\
&= \vc{G}^{m-k}\vc{F}_+(k; k)
+ \vc{G}^{m-1}(\vc{I}-\vc{\Phi}(0))^{-1}\vc{B}(-1)\vc{H}(0;k)
- \vc{u}(m)\vc{\pi}(k),
\label{eq:Hkl}
\end{align}
where the $\vc{u}(m)$ is a vector such that $(\vc{u}(m))_i = \EE[T_0 \mid (X_0,J_0) = (m,i)]$ (see \cite[Lemma~3.5]{Masu21-M/G/1-Subexp}) and
\begin{align}
\vc{u}(m) &
= (\vc{I} - \vc{G}^m)(\vc{I} -\vc{A} -\ol{\vc{m}}_A\vc{g})^{-1}\vc{e}
+ \frac{m}{-\sigma}\vc{e}, \qquad m \in \bbN.
\label{eq:uk}
\end{align}
It also follows from (\ref{eq:Hkl}) that, for any $N \in \bbN$ and $k \in \bbZ_{[0,N]}$,
\begin{align}
&\sum_{n = N + 1}^{\infty} \vc{B}(n)
\left\{ \vc{H}(N; k) - \vc{H}(n; k) \right\}
\nonumber
\\
&\quad=
\sum_{n = N + 1}^{\infty} \vc{B}(n)
\left( \vc{G}^{N-k} - \vc{G}^{n-k} \right)\vc{F}_+(k; k)
\nonumber
\\
&\qquad~~~
+
\sum_{n = N + 1}^{\infty} \vc{B}(n)
\left( \vc{G}^{N-1} - \vc{G}^{n-1} \right)
(\vc{I} - \vc{\varPhi}(0))^{-1}\vc{B}(-1)\vc{H}(0; k)
\nonumber
\\
&\quad\qquad~~
+
\sum_{n = N + 1}^{\infty} \vc{B}(n)
\left\{ \vc{u}(n) - \vc{u}(N) \right\} \vc{\pi}(k).
\label{eq:BH-BH}
\end{align}
Furthermore, \eqref{eq:uk} rewrites the last term in (\ref{eq:BH-BH}) as
\begin{align}
&\sum_{n = N + 1}^{\infty} \vc{B}(n)
\left\{ \vc{u}(n) - \vc{u}(N) \right\}
\nonumber
\\
&\quad= \sum_{n=N+1}^{\infty} \vc{B}(n)(\vc{G}^N - \vc{G}^n)
(\vc{I} -\vc{A} -\ol{\vc{m}}_A\vc{g})^{-1}\vc{e}
+ {1\over{-\sigma}}\sum_{n=N+1}^{\infty} (n - N)\vc{B}(n)\vc{e}
\nonumber
\\
&\quad=
\sum_{n=N+1}^{\infty} \vc{B}(n)(\vc{G}^N - \vc{G}^n)
(\vc{I} -\vc{A} -\ol{\vc{m}}_A\vc{g})^{-1}\vc{e}
+ {1\over{-\sigma}}\ool{\vc{B}}(N - 1)\vc{e},
\label{eq:sumBu-tildeBu}
\end{align}
where the second equality holds due to $\sum_{n=N+1}^{\infty} (n - N)\vc{B}(n)=\ool{\vc{B}}(N - 1)$. Substituting \eqref{eq:sumBu-tildeBu} into \eqref{eq:BH-BH}, and using \eqref{eq:defS}, we obtain
\begin{align}
&\sum_{n = N + 1}^{\infty} \vc{B}(n)
\left\{ \vc{H}(N; k) - \vc{H}(n; k) \right\}
\nonumber
\\
&\quad=
\frac{1}{-\sigma}\ool{\vc{B}}(N - 1)\vc{e}\vc{\pi}(k)
+
\sum_{n = N + 1}^{\infty} \vc{B}(n)
\left( \vc{G}^{N-k} - \vc{G}^{n-k} \right)\vc{F}_+(k; k)
\nonumber
\\
&\qquad+
\sum_{n = N + 1}^{\infty} \vc{B}(n)\left( \vc{G}^{N-1} - \vc{G}^{n-1} \right)
\cdot
(\vc{I} - \vc{\varPhi}(0))^{-1}\vc{B}(-1)\vc{H}(0; k)
\nonumber
\\
&\qquad+
\sum_{n = N + 1}^{\infty} \vc{B}(n)\left( \vc{G}^{N-1} - \vc{G}^{n-1} \right)
\cdot \vc{G}(\vc{I} -\vc{A} -\ol{\vc{m}}_A\vc{g})^{-1}\vc{e}\vc{\pi}(k)
\nonumber
\\
&\quad=
\frac{1}{-\sigma}\ool{\vc{B}}(N - 1)\vc{e}\vc{\pi}(k)
+
\sum_{n = N + 1}^{\infty} \vc{B}(n)
\left( \vc{G}^{N-k} - \vc{G}^{n-k} \right)\vc{F}_+(k; k)
\nonumber
\\
&\qquad+
\sum_{n = N + 1}^{\infty} \vc{B}(n)\left( \vc{G}^{N-1} - \vc{G}^{n-1} \right)
\vc{S}(k),
\qquad k\in\bbZ_{[0,N]}.
\label{eq:BH-BH-2}
\end{align}
Proceeding as in the derivation of (\ref{eq:BH-BH-2}), we have
\begin{align}
&\sum_{n = N + 1}^{\infty} \vc{A}(n)
\left\{\vc{H}(N + \ell; k) - \vc{H}(n + \ell; k)\right\}
\nonumber
\\
&\quad= \frac{1}{-\sigma}\ool{\vc{A}}(N - 1)\vc{e} \vc{\pi}(k)
+
\sum_{n = N + 1}^{\infty} \vc{A}(n)
\left(\vc{G}^{N+\ell-k} - \vc{G}^{n+\ell-k} \right)\vc{F}_+(k;k)
\nonumber
\\
&\qquad+
\sum_{n = N + 1}^{\infty} \vc{A}(n) (\vc{G}^{N +\ell-1}-\vc{G}^{n +\ell-1})\vc{S}(k),\qquad k\in\bbZ_{[0,N]}.
\label{eq:AH-AH-2}
\end{align}
Finally, combining \eqref{eq:piN-pi2} with \eqref{eq:BH-BH-2} and \eqref{eq:AH-AH-2} results in (\ref{eq:piN(k)-pi(k)}). The proof is completed.

\section{Proof of Theorem~\ref{th:subexp}}
\label{proof:th:subexp}

Equation (\ref{eq:subexp-norm-01}) is an immediate consequence of (\ref{eq:subexp-01}), and thus we prove the latter. To do this, we confirm that (\ref{eq:subexp-01}) holds if
\begin{subequations}\label{eqn:210211-02}
\begin{align}
\lim_{N \to \infty}
\vc{\pi}^{(N)}(0)\sum_{n=N+1}^{\infty}
{
\vc{B}(n) ( \vc{G}^{N-k} - \vc{G}^{n-k} )
\over
\ol{F}(N)
}
&=\vc{0},
\label{hm-add-200113-01}
\\
\lim_{N \to \infty}
\sum_{\ell=1}^{\infty}
\vc{\pi}^{(N)}(\ell)
\sum_{n=N+1}^{\infty}
{ \vc{A}(n) (\vc{G}^{N+\ell-k} - \vc{G}^{n+\ell-k})
\over
\ol{F}(N)
}
&=\vc{0},
\label{hm-add-200113-02}
\end{align}
\end{subequations}
for any fixed $k \in \bbZ_+$. It follows from (\ref{eq:piN(k)-pi(k)}) that, for $k\in\bbZ_+$ and $N\in\bbN\cap\bbZ_{\ge k}$,
\begin{align}
&{ \vc{\pi}^{(N)}(k) - \vc{\pi}(k) \over \ol{F}(N) }
\nonumber
\\
&~~= \vc{\pi}^{(N)}(0)
\left[
{1 \over -\sigma}
{\ool{\vc{B}}(N - 1)\vc{e} \over \ol{F}(N) }
\vc{\pi}(k)
+
\sum_{n=N+1}^{\infty}
{
\vc{B}(n) ( \vc{G}^{N-k} - \vc{G}^{n-k} )
\over
\ol{F}(N)
}
\vc{F}_+(k;k)
\right.
\nonumber
\\
&\qquad\qquad\qquad\qquad\quad
\Biggl.
+
\sum_{n=N+1}^{\infty}
{
\vc{B}(n)( \vc{G}^{N-1} - \vc{G}^{n-1} )
\over \ol{F}(N)
}
\vc{S}(k)
\Biggr]
\nonumber
\\
&~~~~
+ \sum_{\ell=1}^{\infty}
\vc{\pi}^{(N)}(\ell)
\left[
{1 \over -\sigma}
{ \ool{\vc{A}}(N - 1)\vc{e} \over \ol{F}(N) }
\vc{\pi}(k)
+ \sum_{n=N+1}^{\infty}
{\vc{A}(n) (\vc{G}^{N + \ell - k} - \vc{G}^{n + \ell - k})
\over
\ol{F}(N)
}
\vc{F}_+(k;k)
\right.
\qquad
\nonumber
\\
&\qquad\qquad\qquad\qquad\quad
\Biggl.
+ \sum_{n=N+1}^{\infty}
{ \vc{A}(n) (\vc{G}^{N + \ell - 1} - \vc{G}^{n + \ell - 1})
\over
\ol{F}(N)
}
\vc{S}(k)
\Biggr].
\label{eq:piN-pi-6}
\end{align}
It also follows from Assumption \ref{assum3} and $F \in \calL$ (see Definition~\ref{defn-long-tailed}) that
\begin{subequations}\label{eq:ABdoublebarP}
\begin{align}
&\lim_{N \to \infty} \frac{\ool{\vc{A}}(N - 1)\vc{e}}{\ol{F}(N)}
=
\lim_{N \to \infty} \frac{\ool{\vc{A}}(N - 1)\vc{e}}{\ol{F}(N-1)}
\frac{\ol{F}(N-1)}{\ol{F}(N)}
= \vc{c}_A,
\label{eq:AdoublebarP}
\\
&\lim_{N \to \infty} \frac{\ool{\vc{B}}(N - 1)\vc{e}}{\ol{F}(N)}
=
\lim_{N \to \infty} \frac{\ool{\vc{B}}(N - 1)\vc{e}}{\ol{F}(N-1)}
\frac{\ol{F}(N-1)}{\ol{F}(N)}
= \vc{c}_B.
\label{eq:BdoublebarP}
\end{align}
\end{subequations}
Using (\ref{eq:ABdoublebarP}) and Theorem~\ref{th:convergence_tvn}, we obtain
\begin{alignat*}{2}
\lim_{N\to\infty}\vc{\pi}^{(N)}(0)
\frac{\ool{\vc{B}}(N - 1)\vc{e}}{-\sigma\ol{F}(N)}\vc{\pi}(k)
&=\frac{\vc{\pi}(0)\vc{c}_B}{-\sigma}\vc{\pi}(k), &\qquad k &\in\bbZ_+,
\\
\lim_{N\to\infty}\sum_{\ell=1}^{\infty}\vc{\pi}^{(N)}(\ell)
\frac{\ool{\vc{A}}(N - 1)\vc{e}}{-\sigma\ol{F}(N)}\vc{\pi}(k)
&=\frac{\ol{\vc{\pi}}(0)\vc{c}_A}{-\sigma}\vc{\pi}(k),&\qquad k &\in\bbZ_+.
\end{alignat*}
Applying these obtained equations and (\ref{eqn:210211-02})  to (\ref{eq:piN-pi-6}), we have (\ref{eq:subexp-01}). Furthermore, the right-hand side of (\ref{eq:subexp-01}) is positive because $\vc{\pi} > \vc{0}$ and either $\vc{c}_A$ or $\vc{c}_B$ is a nonzero vector. Therefore, the proof of (\ref{eq:subexp-01}) is reduced to those of (\ref{hm-add-200113-01}) and (\ref{hm-add-200113-02}).

To complete the proof, we prove (\ref{eqn:210211-02}). It follows from (\ref{eq:ABdoublebarP}) that
\begin{subequations}\label{eq:bar/F}
\begin{align}
&\lim_{N \to \infty} \frac{\ol{\vc{A}}(N){\vc{e}}}{\ol{F}(N)}
=\frac{\ool{\vc{A}}(N)-\ool{\vc{A}}(N-1)}{\ol{F}(N)}
=\vc{0},
\label{eq:AbarP}
\\
&\lim_{N \to \infty} \frac{\ol{\vc{B}}(N)\vc{e}}{\ol{F}(N)}
=\frac{\ool{\vc{B}}(N)-\ool{\vc{B}}(N-1)}{\ol{F}(N)}
=\vc{0}.
\label{eq:BbarP}
\end{align}
\end{subequations}
It also follows from (\ref{defn-ol{A}(k)}), (\ref{eqn:210211-01}), and \eqref{eq:AbarP} that
\begin{align*}
&
\limsup_{N\to\infty}
\left|
\sum_{\ell=1}^{\infty}
\vc{\pi}^{(N)}(\ell)
\sum_{n=N+1}^{\infty}
{ \vc{A}(n)(\vc{G}^{N+\ell-k} - \vc{G}^{n+\ell-k})
\over
\ol{F}(N)
}
\right|
\\
&\quad\leq
\limsup_{N\to\infty}
\sum_{\ell=1}^{\infty}
\vc{\pi}^{(N)}(\ell)
\sum_{n=N+1}^{\infty}
{ \vc{A}(n)\,|\vc{G}^{N+\ell-k} - \vc{G}^{n+\ell-k}|
\over
\ol{F}(N)
}
\\
&\quad\leq
\limsup_{N\to\infty}\vc{e}^{\top}
{
2\ol{\vc{A}}(N)\vc{e}
\over
\ol{F}(N)
}\vc{e}^{\top}
= \vc{0}\quad \mbox{for any fixed $k \in \bbZ_+$},
\end{align*}
which shows that (\ref{hm-add-200113-01}) holds.
Similarly, we have
\begin{align*}
&\limsup_{N\to\infty}
\left|\vc{\pi}^{(N)}(0)
\sum_{n=N+1}^{\infty}
{ \vc{B}(n) (\vc{G}^{N - k} -  \vc{G}^{n - k})
\over
\ol{F}(N)
}
\right|
\\
&\quad\leq
\limsup_{N\to\infty}
\vc{e}^{\top}\frac{2\ol{\vc{B}}(N)\vc{e}}{\ol{F}(N)}\vc{e}^{\top}
= \vc{0} \quad \mbox{for any fixed $k \in \bbZ_+$},
\end{align*}
and thus (\ref{hm-add-200113-02}) holds. The proof is completed.

\section{Subgeometric Functions and Long-tailed Distributions}
\label{sec-subgeo}

This section presents definitions and basic results on the classes of subgeometric functions and long-tailed distributions (see, e.g., \cite{Foss11}). For later use, let $\bbR_+=[0,\infty)$, and let $\ol{F} = 1 - F$ for any distribution function $F$. In addition, we write $f_1(x) = o(f_2(x))$ if $\lim_{x\to\infty}f_1(x)/f_2(x) = 0$.

We first introduce the class of subgeometric rate functions.
\begin{DEF}\label{defn-subgeo}
A function $r:\bbZ_+\to\bbR_+$ is called a {\it subgeometric function} if and only if either
\[
\lim_{k\to\infty}{\log r(k) \over k} = 0 \quad \mbox{or} \quad
\lim_{k\to\infty}{-\log r(k) \over k} = 0.
\]
The class of subgeometric functions is denoted by $\varTheta$.
\end{DEF}

We next introduce some classes of distributions associated with class $\varTheta$.

\begin{DEF}[{}{\cite[Definitions 2.21 and 3.1]{Foss11}}]\label{defn-long-tailed}
\hfill
\begin{enumerate}
\item A distribution function $F$ is said to be {\it long-tailed} if and only if $\ol{F}(k) > 0$ for $k \in \mathbb{Z}_+$ and
\[
\lim_{k \to \infty} \frac{\ol{F}(k + n)}{\ol{F}(k)} = 1, \quad
\forall n \in \mathbb{Z}_+.
\]
\item A distribution function $F$ is said to be {\it subexponential} if and only if $\ol{F}(k) > 0$ for $k \in \mathbb{Z}_+$ and
\[
\lim_{k \to \infty} {\overline{F^{*2}}(k) \over \ol{F}(k)} = 2,
\]
where $F^{*n}$, $n \in \bbN$, denotes the $n$-fold convolution of $F$ itself, i.e.,
\begin{equation*}
F^{*n}(k) =
\begin{cases}
F(k), & k \in \mathbb{Z}_+$, $n=1,
\\
\displaystyle \sum_{\ell = 0}^k F^{*(n - 1)}(k - \ell)F(\ell), & k \in \mathbb{Z}_+$, $n=2,3,\dots.
\end{cases}
%
\end{equation*}
\end{enumerate}
The classes of long-tailed and subexponential distributions are denoted by $\calL$ and $\calS$, respectively.
\end{DEF}

\begin{REM}\label{rem:H=Lambda}
The inclusion relation holds: $\calS \subsetneq \calL$ (see \cite[Lemmas 2.17 and 3.2]{Foss11}). In addition, if $F \in \calL$, then $\ol{F} \in \varTheta$ (see \cite[Proposition~5.8]{Masu21-M/G/1-Subexp}).
\label{rem1}
\end{REM}

Finally, we provide representative examples of subexponential distributions.

\begin{EXA}[{}{Pareto distribution}]\label{exa:Pareto}
The distribution function $F$ such that
\begin{equation*}
\ol{F}(k) = \left({\gamma \over k + \gamma}\right)^\alpha, \qquad \alpha, \gamma \in \bbR_+,
\end{equation*}
is called {\it Pareto distribution}. This distribution has a finite mean if and only if $\alpha > 1$.
\end{EXA}
\begin{EXA}[{}{Heavy-tailed Weibull distribution}]\label{exa:Weibull}
The distribution function $F$ such that
\begin{equation*}
\ol{F}(k) = \re^{-\lambda k^\alpha}, \qquad \lambda \in \bbR_+,\quad 0 < \alpha < 1,
\end{equation*}
is called {\it heavy-tailed Weibull distribution}. This distribution always has a finite mean.
\end{EXA}
%




\section*{Acknowledgments}
The research of the second author was supported in part by JSPS KAKENHI Grant Number JP21K11770.

%
%
%
\bibliographystyle{plain} 


\end{document}